

\input amssym.def

 \input amssym.tex
 \font\newrm =cmr10 at 24pt
\def\bul{\raise .9pt\hbox{\newrm .\kern-.105em } }

 \def\fr{\frak}

 \baselineskip 20pt
 \magnification=1200
 \def\h{\hbox{ }}
 
 \def\p{{\fr p}}

 \def\ss{{\fr s}}
 \def\k{{\fr k}}

 \def\g{{\fr g}}
 \def\v{{\fr v}}

 \def\<{\le}
 \def\>{\ge}

 \def\s{{\h\subset\h}}

 \def\mapright#1
  {\smash{\mathop
  {\longrightarrow}
  \limits^{#1}}}

 \def\kk#1{{\kern .4 em} #1}


\font\twelverm=cmr12 at 14pt
\font\authorfont=cmr9
\font\ninerm=cmr9
\centerline{\twelverm The Weyl algebra and the structure of all }
 \centerline {\twelverm Lie superalgebras of Riemannian type} 
\vskip 1.5pc
\baselineskip=11pt
\vskip8pt
\centerline{\authorfont BERTRAM KOSTANT\footnote*{\ninerm
Research supported in part by NSF grant
DMS-9625941 and in part by the \hfil\break KG\&G Foundation}}\vskip 2pc
\baselineskip 20pt
\centerline{\bf 0. Introduction}\vskip 1.5pc
0.1. All vector spaces are over $\Bbb C$ and are assumed to be finite dimensional unless
manifestly otherwise (e.g., tensor algebras, Weyl algebras, symmetric algebras, universal
enveloping algebras of Lie algebras). A bilinear form 
$(u,v)$ on a $\Bbb Z_2$-graded vector space $V = V_0 + V_1$ 
will be called here supersymmetric in case $(x,y) =
(-1)^{|x||y|} (y,x)$ for homogeneous elements $x,y\in V$ and the even component
$V_0$ is orthogonal to the odd component $V_1$. In case $V=\g$ is a Lie superalgebra
with $\g_0 = V_0$ and $V_1=\g_1$, such a bilinear form, now denoted by $B_{\g}$, will be
called super-Riemannian in case
$B_{\g}$ is non-singular and
$ad$-invariant. In such a case
$B_{\g_0} = B_{\g}|\g_0$ is non-singular, symmetric (in the usual sense) and
$ad$-invariant for the ordinary Lie algebra $\g_0$. Also $B_{\g_1} = B_{\g}|\g_1$ is
non-singular and alternating in the usual sense. Furthermore if $Sp(\g_1)\s Aut\, \g_1$ is
the symplectic group with respect to $B_{\g_1}$, then one has a representation
 $$\nu:\g_0\to
Lie\,Sp(\g_1)\eqno (0.1)$$ where for $x\in \g_0,\, y\in \g_1$, one has $\nu(x)y = [x,y]$. 
A Lie superalgebra which admits a super-Riemannian bilinear form will be said to be of
Riemannian type. \vskip 1pc {\bf Examples.}\item {[1]} In his classification of all
simple Lie superalgebras Victor Kac (see [Kac]) isolated a subfamily which he referred to
as basic classical Lie superalgebras. These are all of Riemannian type. Besides ordinary
simple Lie algebras they list as $A(m,n)$, $B(m,n)$, $C(n)$, $D(m,n)$, $D(2,1,\alpha)$,
$F(4),$ and
$G(3)$.

\item{[2]} If $\ss$ is any Lie superalgebra and $\ss^*$ is its dual we may form the
semi-direct sum $\g= \ss \oplus \ss^*$ with $\ss$ operating on $\ss^*$ by the coadjoint
representation and with $\ss^*$ taken to be abelian. Then $\g$ has the structure of a Lie
superalgebra of Riemannian type.\vskip 1pc 0.2. Let $\v$ be a vector
space with a non-singular alternating bilinear form $B_{\v}$ and let $Sp(\v)$ be the
symplectic group of
$\v$ with respect to $B_{\v}$. Let $\g_0$ be an (ordinary) Lie algebra
which admits a non-singular symmetric bilinear form $B_{\g_0}$. Assume one has 
a representation $$\nu:\g_0\to Lie\,Sp(\v)\eqno (0.2)$$ As noted in (0.1) such a
representation arises from a Lie superalgebra of Riemannian type. The main theorem
(Theorem 0.1 below) gives a necessary and sufficient condition, purely in terms of
the symplectic representation $\nu$, that (0.2) be of the form (0.1). 

Let $W(\v)$ be the Weyl algebra over $\v$ (with respect to $B_{\v}$) so that $W(\v)$ is the
free algebra over
$\v$ subject to the relations $u\,v- v\,u = 2(u,v)$ for all $u,v\in \v$. Let $S(\v)=
\sum_{n=0}^{\infty}S^n(\v)$ be the symmetric algebra over $\v$. The bilinear form
$B_{\v}$ on $\v$ extends in a natural way to a non-singular bilinear form $B_{S(\v)}$
on $S(\v)$. There is a natural identification of the underlying spaces of
$W(\v)$ and
$S(\v)$ so that we will understand that $W(\v)$ has two multiplicative structures. The
identification is such that
$u^n$ for $u\in \v$ and $n\in \Bbb Z_+$ has the same meaning in both structures. The Weyl
product of 
$w,z\in W(\v)$ is denoted by $w\,z$ and the product in $S(\v)$ will be denoted by
$w\cdot z$. Commutation using Weyl multiplication defines a Lie algebra structure on
$S^2(\v)$ and a Lie algebra isomorphism $$S^2(\v)\to Lie\,Sp(\v)\eqno (0.3)$$ The
representation $\nu$ then uniquely lifts to a Lie algebra homomorphism $$\nu_*:\g_0\to
S^2(\v)\eqno (0.4)$$ so that $\nu_*$ composed with the isomorphism (0.3) is $\nu$. With
respect to the bilinear forms $B_{\g_0}$ and $B_{S(\v)}|S^2(\v)$ the map $\nu_*$ admits a
transpose $$\nu_*^t:S^2(\v)\to \g_0\eqno (0.5)$$  Now let $U(\g_0)$ be the universal
enveloping algebra of $\g_0$. The Lie algebra homomorphism (0.4) then extends uniquely to an
algebra homomorphism $$\nu_*:U(\g_0)\to W(\v)\eqno (0.6)$$ Let $x_i,\,i=1,\ldots,k,$ be an
orthonormal basis of $\g_0$ with respect to $B_{\g_0}$. The Casimir element
$Cas_{\g_0}\in Cent\,U(\g_0)$ corresponding to $B_{\g_0}$ is given by $Cas_{\g_0} =
\sum_{i=0}^k\,x_i^2$. If we identify $S^0(\v)$ with $\Bbb C$ one readily shows that
$$\nu_*(Cas_{\g_0}) \in S^4(\v) + \Bbb C\eqno (0.7)$$ 

Now let $\g = \g_0 \oplus \v$. Regarding $\g$ as a $\Bbb Z_2$-graded vector space with
$\g_0$ and $\v$ respectively as the even and odd components, let $B_{\g}$ be the
non-singular supersymmetric bilinear form on $\g$ defined so that $B_{\g}\vert \g_0 =
B_{\g_0}$ and $B_{\g}\vert \v = B_{\v}$. We will now say that the pair $(\nu,B_{\g})$ is
of super Lie type in case there is a Lie superalgebra structure on $\g$ so that (1)
$\g_0$, as a Lie algebra, is the even component of $\g$, (2) $\v= \g_1$, and for any
$x\in \g_0,\,y\in \v$ one has $[x,y] = \nu(x)y$, and (3) $B_{\g}$ is a super-Riemannian
bilinear form on $\g$ (so that in particular $\g$ is of Riemannian type). The
following main theorem is proved as Theorem 2.8 in \S 2. \vskip 1pc {\bf Theorem 0.1.} {\it
Let the pair
$(\nu,B_{\g})$, as above, be arbitrary. Then
$(\nu,B_{\g})$ is of super Lie type if and only if $\nu_*(Cas_{\g_0})$ is a constant.
That is, in the notation of (0.7),  $(\nu,B_{\g})$ is of super Lie type if and only if
the component of $\nu_*(Cas_{\g_0})$ in $S^4(\v)$ vanishes. Furthermore in such a case
the Lie superalgebra structure on $\g$ is unique. In fact for $y,y'\in \v$ the element
$[y,y']\in \g_0$ is given by $$[y,y'] = 2\,\nu_*^t(y\cdot y')\eqno (0.8)$$} \vskip 1pc
{\bf Remark 0.2.} We may restate Theorem 0.1 in the language of representation theory.
Recall that $\v$ is even dimensional, say $dim\,\v = 2m$. Then if 
$Diff_{poly}\Bbb C^m$ is the algebra of all differential operators on $\Bbb C^m$ with
polynomial coefficients one has an algebra isomorphism $$W(\v)\to Diff_{poly}\Bbb
C^m\eqno (0.9)$$ If we compose (0.9) with (0.6) one has a homomorphism $$U(\g_0)\to
Diff_{poly}\Bbb C^m\eqno (0.10)$$ The statement of Theorem 0.1 is that $(\nu,B_{\g})$ is
of super Lie type if and only if the image of the Casimir element $Cas_{\g_0}$ in (0.10) is
just multiplication by a constant. \vskip 1pc 0.3. The statement of Theorem 0.1 was
suggested to us by an analogous statement for ordinary Lie algebras which we proved some
time ago. If $\g$ is an ordinary Lie algebra with a non-singular $ad$-invariant symmetric
bilinear form $B_{\g}$ and $\k\s \g$ is a Lie subalgebra such that $B_{\k} = 
B_{\g}\vert \k$, is non-singular let $\p$ be the $B_{\g}$ orthocomplement of $\k$ in $\g$
and let $B_{\p} = B_{\g}\vert \p$. Let $SO(\p)$ be the orthogonal group with respect to
$B_{\p}$. The relation $[\k,\p]\s \p$ defines a representation $$\nu:\k\to
Lie\,SO(\p)\eqno (0.11)$$ One question we considered and solved in the early 1960's was to
characterize the case where $\k$ is a symmetric subalgebra of $\g$ (i.e., the case where
$[\p,\p]\s \k$) purely in terms of the representation $\nu$. The characterization is the
analogue of Theorem 0.1 where the Clifford algebra $Cliff(\p)$ replaces the Weyl algebra
$W(\v)$ and the exterior algebra
$\wedge\p$ replaces the symmetric algebra $S(\v)$. The condition is again that
$\nu_*(Cas_{\k})$ be a constant. For the statement see p.152 in [Con]. For the proof see
Theorem 1.59 in [Kos]. \vskip 1pc {\bf Remark 0.3.} The constancy of 
$\nu_*(Cas_{\k})$ in case rank $\k$ = rank $\g$ was used by Parthasarathy to construct
discrete series representations. See Lemma 2.2
in [Parth]. One is left to wonder whether Theorem 0.1 will have
analogous applications for representations of Lie superalgebras of Riemannian type.
\vskip 1.5pc 


\hsize = 31pc
\vsize = 45pc
\overfullrule = 0pt
	\centerline{\bf 1. Lie superalgebras }\vskip 1.5pc
1.1. If $$V = V_0 + V_1$$ is a
$\Bbb Z_2$-graded vector space and $v\in V$, 
then the use of the notation $|v|$ implicitly assumes
that $v$ is homogeneous and $|v| = 0$ or $1$
according as $v$ is even, i.e., $v\in V_0$ or odd,
i.e., $v\in V_1$. If $i,j\in \{0,1\}$ represent homogeneity indices in $\Bbb Z_2$-graded
vector spaces then $i+j$ is taken in $\Bbb Z_2$. A 
$\Bbb Z_2$-graded vector space
$\g =
\g_0 +
\g_1$ is a Lie superalgebra if the there is a bilinear operation $[x,y]$ on
$\g$, where $$[\g_i,\g_j]\s \g_{i+j}\eqno (1.1)$$ such that for $x,y,z\in\g$ 
$$\eqalign{(A)&\,\,\,\,[x,y] = -(-1)^{|x||y|}[y,x]\,\,\hbox{and}\cr
(B)&\,\,\,\,(-1)^{|x||z|}[x,[y,z]] +  (-1)^{|z||y|}[z,[x,y]] + (-1)^{|y||x|}[y,[z,x]] =
0\cr}$$  One readily shows that $(A)$ and $(B)$ are equivalent to $(A)$ and $(B')$
where $(B')$ is the super derivation condition $$[x,[y,z]] = [[x,y],z] +
(-1)^{|x||y|}[y,[x,z]]\eqno (1.2)$$ If
$V$ is a
$\Bbb Z_2$-graded vector space then
$End\,V$ has the structure of a  Lie superalgebra where for $a\in End\,V$ and
$v\in V$ one has $|av| = |a| + |v|$ and for $a,b\in End\,V $ one has $[a,b] =
ab-(-1)^{|a||b|}ba$. If $\g$ is Lie superalgebra then a super representation of $\g$ on $V$
is a homomorphism $\pi:\g\to End\,V$ of Lie superalgebras which preserves the grading.

Let $\g$ be a Lie superalgebra and consider the linear map $$ad:\g\to End\,\g,\,\hbox{where
for
$x,y\in\g$},\,\,ad\,x(y) = [x,y]$$ Then one readily verifies \vskip 1pc {\bf Remark 1.1.} The
map $ad$ is a super representation of $\g$ on itself (referred to as the adjoint
representation). It follows from (1.1) and (1.2) that the restriction  of the bracket operation
to $\g_0$ defines an ordinary Lie algebra structure on $\g_0$ and for $x\in \g_0$ the map
$$x\mapsto ad\,x\vert
\g_1\,\,\hbox{defines an ordinary representation of $\g_0$ on $\g_1$}\eqno (1.3)$$\vskip 1pc
1.2. Let $V$ be a $\Bbb Z_2$-graded vector space. The value taken on an ordered pair $u,v\in
V$ by a bilinear form $B$ on $V$ will be denoted by $(u,v)$ if there is no danger of
confusion. A bilinear form $B$ on $V$ will be said to be supersymmetric if
$$\eqalign{(\alpha)&\,\,(u,v) = (-1)^{|u||v|}(v,u),\,\,\hbox{for all homogeneous $u,v\in
V$}\cr (\beta)&\,\,(V_i,V_j) = 0\,\,\hbox{if}\,\,i\neq j\cr}\eqno (1.4)$$ If $\pi$ is a
representation of a Lie superalgebra $\g$ on $V$ then we will say that $B$ is invariant
under $\pi$ (or $\g$ if
$\pi$ is understood) in case $$(\pi(x)u,v) + (-1)^{|x||u|}(u,\pi(x)v) = 0\eqno (1.5)$$ for all
$v\in V,\,\,
\hbox{and homogeneous}\,\,x\in \g, u\in V$.

Let $V$ be a $\Bbb Z_2$-graded vector space. The dual space $V^*$ to $V$ has an obvious
$\Bbb Z_2$-grading. For $v\in V$ and $w\in V^*$ we will denote the pairing of $w$ and $v$
into $\Bbb C$ by $\langle w,v\rangle$. If $a\in End\,V$ then $a^t\in End V^*$ is defined so
that $\langle a^tw,v\rangle = \langle w,av\rangle$ for any $v\in V$ and any $w\in V^*$.
One readily shows that if
$\pi$ is a super representation of a Lie superalgebra on $V$, then $\pi^*$ is a super
representation of $\g$ on $V^*$ (referred to as the representation contragredient to
$\pi$) where for homogeneous $x\in \g$ and homogeneous $w\in V^*$ one has $$\pi^*(x) w =
-(-1)^{|x||w|}\pi(x)^t w\eqno (1.6)$$ \vskip .5pc {\bf Remark 1.2}. Note that it is immediate
from (1.5) and (1.6) that if $V$ is finite dimensional and $B$ is a non-singular
$\pi$-invariant supersymmetric bilinear form on $V$ then the linear isomorphism $\xi:V\to
V^*$ defines a representation equivalence of $\pi$ and $\pi^*$, where for any $u,v\in V$ one
has $\langle \xi u,v\rangle = (u,v)$. \vskip 1pc
1.3. Let $V$ be a finite dimensional $\Bbb Z_2$-graded vector space. The supertrace on
$End\,V$ is a linear map $$str:End\,V\to \Bbb C\eqno (1.7)$$ defined so that $str\,a = 0$ if
$a\in (End\,V)_1$ and for $a\in (End\,V)_0$ one has $$str\,a = tr\,a|V_0 -tr\,a|V_1$$
Regarding $\Bbb C$ as a 1-dimensional Lie algebra (and hence a fortiori a Lie
superalgebra) one readily shows that (1.7) is a homomorphism of Lie superalgebras. In
particular for $a,b\in End\,V$, $$\eqalign{str\,[a,b] &= str\,ab - (-1)^{|a||b|}ba\cr 
&=0\cr}\eqno (1.8) $$ But for $a,b,c\in End\,V$ one has $$[a,bc] = [a,b]c +
(-1)^{|a||b|}b[a,c]\eqno (1.9)$$ Let $B_{End\,V}$ be the bilinear form on $End\,V$ defined
by putting $(b,c) = str\,bc$. It is trivial to see that $B_{End\,V}$ is non-singular. As a
consequence of (1.8) and (1.9) one then recovers the well-known
\vskip 1pc {\bf Proposition 1.3.} {\it Let $V$ be a finite dimensional $\Bbb Z_2$-graded
vector space.  Then $B_{End\,V}$ is a non-singular supersymmetric bilinear form on the Lie
superalgebra $End\,V$. Furthermore $B_{End\,V}$ is invariant under the (super) adjoint
representation of $End\,V$ on itself.} \vskip 1pc In the terminology of the introduction,
which shall be retained from now on, $End\,V$ (wth the bracket structure defined above) is a Lie
superalgebra of Riemannian type and $B_{End\,V}$ is a super-Riemannian bilinear form on
$End\,V$.  But now if
$\pi:\g\to End\,V$ is any super representation of a Lie superalgebra on a $\Bbb
Z_2$-graded vector space $V$ then $B_{End\,V}$, by Proposition 1.3, induces a supersymmetric
bilinear form
$B_{\pi}$ on $\g$ given by $$(x,y) = str\,\pi(x)\pi(y)\eqno (1.10)$$ which is invariant
under the adjoint representation of $\g$. If
$\pi$ is the adjoint representation of
$\g$, then one will refer to $B_{ad}$ as the super Killing form. 

Now if $\g$ is any Lie superalgebra let $$\nu:\g_0\to End\,\g_1\eqno (1.11)$$ be the
representation of the ordinary Lie algebra $\g_0$ on the vector space $\g_1$ given by
$\nu(x)y = [x,y]$ for $x\in \g_0$ and $y\in \g_1$. In this paper the main result will
give the structure of an arbitrary Lie superalgebra $\g$ of Riemannian type. The structure
theorem will depend only on the ordinary Lie algebra $\g_0$ and a special property of a
vector space representation $\nu$ of $\g_0$. The following are a few examples of Lie
superalgebras of Riemannian type. \vskip 1pc {\bf Examples.}

 \item {(1)} The Kac list of basic classical Lie superalgebras (see Introduction, \S 0.1). One
notes that, for the case (i.e., $A_{m,n}$) of Proposition 1.1, if
$dim\,V_0=dim\,V_1$ then even though $B_{End\,V}$ is non-singular, the super Killing form is
identically 0. 

\item {(2)} Let
$B$ be any adjoint invariant, supersymmetric bilinear form on a Lie
superalgebra $\g$. For example take $B = B_{\pi}$ in the notation of (1.10). Let $rad\,B$ be
the radical of $B$ so that $rad\,B$ is a Lie superideal of $\g$. Then $\g/rad\,B$ is a Lie
superalgebra of Riemannian type and $B'$ is a super-Riemannian bilinear form on $\g/rad\,B$
where $B'$ is the bilinear form on $\g/rad\,B$ induced by $B$. One notes
that if $\g$ is simple either $rad\,B = 0$ or $rad\,B =\g$. 

\item {(3)} Let $\ss$ be any finite dimensional Lie superalgebra. The coadjoint
representation, $coad$, of $\ss$ is the representation of $\ss$ on the dual space
$\ss^*$ to $\ss$, which is contragredient to $ad$. Let $\g = \ss \oplus \ss^*$. One defines
a Lie superalgebra structure on $\g$ so that $\ss$ is a subalgebra, $[\ss^*,\ss^*]=0$, and
$[x,y] = coad\,x(y)$ for $x\in \ss$ and $y\in \ss^*$. In such a case put 
$[y,x] = -(-1)^{|x||y|}[x,y]$. One defines a bilinear form
$B_{\g}$ on $\g$ so that $\ss$ and $\ss^*$ are isotropic subspaces and $(y,x) = \langle
y,x\rangle$ for $y\in \ss^*$ and $x\in \ss$. In such a case one puts $(x,y) =
(-1)^{|x||y|}(y,x)$. It is readily verified that $B_{\g}$ is a super-Riemannian bilinear form on
$\g$. \vskip 1.5pc \centerline{\bf 2. The Weyl algebra and the main theorem}\vskip 1.5pc 2.1. Let
$\v$ be an even dimensional vector space and let $B_{\v}$ be a non-singular alternating
bilinear form
$(u,v)$ on $\v$. Let $T(\v)$ be the tensor algebra over $\v$. Let $I_W(\v)$ (resp. $I_S(\v)$) be
the ideal in $T(\v)$ generated by all elements of the form $u\otimes v-v\otimes u - 2(u,v)$ (resp.
$u\otimes v-v\otimes u $) for $u,v\in \v$. The Weyl algebra
$W(\v)= T(\v)/I_W(\v)$ (resp. the symmetric algebra $S(\v) = T(\v)/I_S(\v)$) and one has a
quotient homomorphism $\sigma_W:T(\v)\to W(\v)$ (resp. $\sigma_S:T(\v)\to S(\v)$). We may 
identify $\v$ with $\sigma_W(\v)$ and with $\sigma_S(\v)$ so that $S(\v)$ is the commutative
free algebra generated by $\v$ and $W_{\v}$ is the algebra generated by $\v$ where one has the
relation $$uv-vu = 2(u,v)\eqno (2.1)$$ for $u,v\in \v$. If
$D_T$ is a derivation of $T(\v)$ which stabilizes both $I_W(\v)$ and $I_S(\v)$ then $D_T$
descends to a derivation $D_W$ of $W(\v)$ and a derivation $D_S$ of $S(\v)$. We consider two
examples. Let $Sp(\v)\s Aut\,\v$ be the Lie algebra of the symplectic group with respect
to $B_{\v}$. If $\alpha\in Lie\,Sp(\v)$ then $\alpha$ extends uniquely to a derivation 
$\tau_T(\alpha)$ of $T(\v)$. The extension clearly stabilizes both $I_W(\v)$ and $I_S(\v_)$
so that $\alpha$ uniquely extends to the derivations $\tau_W(\alpha)$ and
$\tau_S(\alpha)$ of $W(\v)$ and $S(\v)$, respectively. The symmetric algebra $S(\v)$ is of
course graded and $\tau_S(\alpha)$ clearly stabilizes the homogeneous components
$S^k(\v),\,k=0,1,...,$. For the second example one notes that if $u\in \v$, the unique derivation
$\iota_T(u)$ (of degreee $-1$) of $T(\v)$ such that $\iota_T(u)v= (u,v)$ stabilizes both
$I_W(\v)$ and $I_S(\v_)$ and hence descends to derivations $\iota_W(u)$ and $\iota_S(\v)$ of
$W(\v)$ and $S(\v)$, respectively.

One naturally extends the bilinear form $B_{\v}$ on $\v$ to a non-singular bilinear form
$B_{S(\v)}$ on $S(\v)$ so that $S^m(\v)$ is orthogonal to $S^n(\v)$ for $m\neq n$ and if
$u_i,v_j\in \v$ for $i,j=1,\ldots,n$, then $$(u_1\cdots u_n,v_1\cdots v_n) =
\sum_{\sigma} (u_1,v_{\sigma 1})\cdots  (u_n,v_{\sigma n})\eqno (2.2)$$ where the sum is over
the permutation group of $\{1,\ldots,n\}$.  \vskip 1pc {\bf Remark 2.1.} Note that $B_{S(\v)}$
is symmetric on $S^n(\v)$ if $n$ is even and alternating on $S^n(\v)$ if $n$ is odd.\vskip 1pc 

For each $u\in
\v$ let $\varepsilon_S(u)$ be the linear operator on $S(\v)$ of multiplication by $u$. It is
immediate from (2.2) that for any $u\in \v$ and $w,z\in S(\v)$ that $$(\varepsilon_S(u)w,z)
= (w,\iota_S(u)z)\eqno (2.3)$$ Also for $u,v\in \v$ one readily has that the commutation
relations
$$\eqalign{[i_S(u),\varepsilon_S(v)]& = [\varepsilon_S(u),\iota_S(v)]\cr
&= (u,v)\cr}\eqno (2.4)$$ It follows therefore that if $\gamma(u) = \varepsilon_S(u) +
\iota_S(u)$, then for $u,v\in \v$ one has $$[\gamma(u),\gamma(v)] = 2(u,v)\eqno (2.5)$$
The linear map $\v\to End\,S(\v)$ where $u\mapsto \gamma(u)$ naturally extends to a
homomorphism
$T(\v)\to End\,S(\v)$ which, by (2.5), descends to a homomorphism $$\gamma:W(\v)\to
End\,S(\v)\eqno (2.6)$$ defining the structure of a $W(\v)$-module on $S(\v)$. Now let
$\eta:W(\v)\to S(\v)$ be the linear map defined by putting $$\eta(w) = \gamma(w) 1\eqno
(2.7)$$ Let $\Sigma(\v)\s T(\v)$ be the space of symmetric tensors. Clearly
$\sigma_W:\Sigma(\v)\to W(\v)$ and $\sigma_S:\Sigma(\v)\to S(\v)$ are linear isomorphisms. If
$u\in
\v$ then since
$(u,u) = 0$, it follows easily that
$\eta$ is a linear isomorphism since one clearly has,
for any $n\in \Bbb Z_+$,
$$\eta(u^n) = u^{[n]}\eqno (2.8)$$ where, to avoid confusion, we have written $u^n$ and
$u^{[n]}$, respectively, for the $n^{th}$ power of $u$ in $W(\v)$ and $S(\v)$.  Henceforth we
will identify $W(\v)$ with
$S(\v )$, using
$\eta$, and understand that there are two multiplications on $W(\v)$, Weyl multplication, $wz$,
and the commutative multiplication of $S(\v)$ denoted now by $w\cdot z$. \vskip 1pc {\bf
Proposition 2.2.} {\it Let $u\in \v$. Then $\iota_W(u) = \iota_S(u)$ so that if we put
$\iota(u)=\iota_W(u) = \iota_S(u)$ then $\iota(u)$ is a derivation of both the Weyl and
symmetric algebra structures on $W(\v)$. 

Let $\alpha\in Lie\,Sp(\v)$. Then $\tau_W(\alpha)=\tau_S(\alpha)$ so that if we put
$\tau(\alpha) = \tau_W(\alpha)=\tau_S(\alpha)$ then $\tau(\alpha)$ is a derivation of both
the Weyl and
symmetric algebra structures on $W(\v)$.} \vskip 1pc {\bf Proof.} It follows from (2.8) that
$$\eta\, \sigma_W = \sigma_S\eqno (2.9)$$ on the subspace $\Sigma(\v)\s T(\v)$ of
symmetric tensors. But then the proposition is a consequence of (2.9) since
$\iota_T(u)$ and $\tau_T(\alpha)$ stabilize $\Sigma(\v)$. Indeed if $D_T,D_W$ and $D_S$ are
as in the beginning of \S 2.1 then $D_W\, \sigma_W = \sigma_W\, D_T$ and $D_S\,\sigma_S =
\sigma_S\, D_T$ on $T(\v)$. But if $D_T$ stabilizes $\Sigma(\v)$ then (2.9) implies $$\eta\,
\sigma_W\, D_T=
\sigma_S\, D_T$$ on $\Sigma(\v)$. Hence on $\Sigma(\v)$, $$\eqalign{\eta\, D_W\,
\sigma_W &= D_S\, \sigma_S\cr &= D_S\, \eta\, \sigma_W\cr}$$ Thus $$\eta\, D_W = D_S\, \eta\eqno
(2.10)$$ on
$W(\v)$. QED\vskip 1pc 2.2. For any element $w\in S^2(\v)$ let $ad\,w$ be the operator on
$W(\v)$ defined so that $ad\,w(z) = [w,z]$ where $z\in W(\v)$ and commutation is with
respect to Weyl multiplication. \vskip 1pc {\bf Lemma 2.3.} {\it Let $w\in S^2(\v)$. Then
$\v = S^1(\v)$ is stable under $ad\,w$ and if $v\in \v,w\in S^2(\v)$ one has
$$ad\,w(v) = -2\,\iota(v)w\eqno (2.11)$$ Furthermore $ad\,w\vert \v \in Lie\,Sp(\v)$ and
if $u,v\in \v$ one has $$(u,ad\,w(v)) = -2\iota(u)\iota(v)w\eqno (2.12) $$}\vskip 1pc {\bf
Proof.} Since
$S^2(\v)$ is clearly spanned by elements of the form $t^2$ for $t\in \v$ to prove (2.11) it
suffices to take
$w=t^2$. But then for
$v\in \v$ $$\eqalign{-2\iota(v)t^2&= -4(v,t)\,t\cr &=4(t,v)\,t\cr
&= t^2\,v-v\,t^2\cr}$$ This establishes (2.11). Now let $u,v\in \v$. But then (2.12)
follows from (2.11). However since $\iota(u)$ and $\iota(v)$ clearly commute (this follows
for example from (2.3)) one has $(u,ad\,w(v)) = (v,ad\,w(u))$. But of course
$(v,ad\,w(u))= - (ad\,w(u),v)$ so that $ad\,w\vert \v \in Lie\,Sp(\v)$. QED \vskip 1pc  Now
let $$A:S^2(\v)\to Lie\,Sp(\v)\eqno (2.13)$$ be the map defined by putting $A(w)=
(ad\,w)\vert \v$. \vskip 1pc {\bf Proposition 2.4.} {\it The subspace $S^2(\v)\s W(\v)$ is a
Lie algebra under Weyl commutation. Moreover $A$ is a Lie algebra isomorphism. In addition if
$w\in S^2(\v)$ and $\alpha = A(w)$, then $$\tau(\alpha) = ad\,w\eqno (2.14)$$} \vskip 1pc {\bf
Proof.} Let $w\in S^2(\v)$ and let $\alpha = A(w)\in Lie\,Sp(\v)$. Obviously $ad\,w$ is a
derivation of $W(\v)$. But then $ad\,w = \tau_W(\alpha)$ by definition of $\tau_W$. Hence 
$ad\,w =\tau(\alpha)$ by Proposition 2.2. Consequently  $ad\,w = \tau_S(\alpha)$ so that
$S^n(\v)$ is stable under $ad\,w$ for all $n\in \Bbb Z_+$. In particular this is the case for
$n=2$ so that $S^2(\v)$ is a Lie subalgebra of $W(\v)$. The Jacobi identity for 
commutation in $W(\v)$ then implies that $A$ is a Lie algebra homomorphism. The map $A$ is
injective by (2.12). It is then surjective by dimension. QED\vskip 1pc 
Identify $S^0(\v)$ with $\Bbb C$. For any $w\in W(\v)$ let $\chi(w)\in \Bbb C$ be the component
in $S^0(\v)$ with respect to the graded structure on $W(\v)$ defined by the homogeneous spaces
$S^n(\v),\,\,n\in \Bbb Z_+$. \vskip 1pc {\bf Proposition 2.5.} {\it The bilinear form
$B_{S(\v)}$ on $W(\g)$ may be given by $$(y,z) = \chi(y\,z)\eqno (2.15)$$ for any $y,z\in
W(\v)$. Furthermore $B_{S(\v)}$ is invariant under $ad\,w$ for any $w\in S^2(\v)$.} \vskip 1pc
{\bf Proof.} Let $y,z\in W(\v)$. By (2.7) $$\eqalign{y\,z& = \gamma(yz)1\cr
&=\gamma(y)\gamma(z)1\cr &= \gamma(y)z\cr}\eqno (2.16)$$ By linearity, to prove (2.15) it
suffices to take $y= u^n,\,z=v^m$ for
$u,v\in
\v$ and $n,m\in \Bbb Z$. But by (2.16) $$y\,z =
(\varepsilon_S(u) +
\iota_S(v))^n\,v^m\eqno (2.17)$$ But then it is immediate from (2.17) that $\chi(y\,z) = 0$ if
$n\neq m$ and if $n=m$, (2.17) implies $$\eqalign{\chi(y\,z) &= (\iota_S(u))^m v^m\cr
&= m!(u,v)^m\cr
&= (y,z)\cr}$$ This establishes (2.15). But now $B_{S(\v)}$ is invariant under $ad\,w$ for
$w\in S^2(\v)$ by (2.15) since $ad\,w$ derivates Weyl multiplication and, also, by (2.14),
$ad\,w$ preserves the gradation of $S(\v)$. QED\vskip 1pc

Clearly there
exists a unique anti-automorphism $J_T$ of the tensor algebra $T(\v)$ such that $J_T(u) = iu$
for any $u\in \v$. One notes that $J_T$ fixes any element of the form $u\otimes v - v\otimes u -
2(u,v)$ where $u,v\in \v$ so that $J_T$ stabilizes the ideal $I_W(\v)$. Consequently $J_T$
descends to an anti-automorphism $J$ of order 4 on the Weyl algebra $W(\v)$ and equals the
scalar $i$ on $\v$. Clearly then for any $n\in \Bbb Z_+$, $$J =
i^n\,\,\hbox{on}\,\,S^n(\v)\eqno (2.18)$$ since $S^n(\v)$ is spanned by all elements of the form
$u^n$ for $u\in \v$. One notes then that
$J^2$ is an automorphism of order 2 on $W(\v)$ and $$J^2 = (-1)^n 
\,\,\hbox{on}\,\,S^n(\v)\eqno (2.19)$$ Let $W^{even}(\v) = \sum_{n=0}^{\infty} S^{2n}(\v)$.
Identify $S^{0}(\v)$ with $\Bbb C$. \vskip 1pc {\bf Proposition 2.6.} {\it $W^{even}$ is the
(associative) subalgebra of $W(\v)$ generated by $S^2(\v)$. Furthermore for any $w\in S^2(\v)$
one has $$w^2 \in S^4(\v) + \Bbb C\eqno (2.20)$$} \vskip 1pc {\bf Proof.} It is clear from
(2.19) that $W^{even}(\v)$ is the space of fixed elements for the automorphism $J^2$. Hence
$W^{even}(\v)$ is a subalgebra of $W(\v)$. However clearly $W^{even}(\v)$ is spanned by all
elements of the form $u^{2n}= (u^2)^n$ for all $u\in \v$. Hence $W^{even}(\v)$ is generated by
$S^2(\v)$. But now it follows immediately from (2.16) and the formula for $\gamma(v)$ for
$v\in \v$ that $$(S^2(\v))^2 \s S^4(\v) + S^2(\v) + \Bbb  C\eqno (2.21)$$ We have only to
show that if $w\in S^2(\v)$, the component of $w^2$ in $S^2(\v)$ on the right side of (2.21)
vanishes. But $J =-1$ on $S^2(\v)$ and $J=1$ on $S^4(\v) + \Bbb  C$. But $J(w^2) = w^2$. Thus
the component of $w^2$ in $S^2(\v)$ vanishes. QED\vskip 1pc 2.3. Now assume that $\g_0$ is
some (ordinary) Lie algebra with a nonsingular $ad$-invariant symmetric
bilinear form $(x,y)$ denoted by $B_{\g_0}$. Assume also that $\v$ is an even dimensional
vector space with a non-singular alternating bilinear form $B_{\v}$ and that $$\nu:\g_0\to
Lie\,Sp(\v)$$ is a symplectic representation (with respect to $B_{\v}$) of $\g_0$ on $\v$.
\vskip 1pc {\bf Remark 2.7.} Note that one is presented with these assumptions in case $\g =
\g_0 + \g_1$ is a Lie superalgebra of Riemannian type. In such a case let $B_{\g}$ be a
super Riemannian bilinear form on $\g$. Indeed let $B_{\g_0} = B_{\g}\vert \g_0$ and
$B_{\g_1}=  B_{\g}\vert \g_1$. Then the assumptions are satisfied where $\v = \g_1$, $B_{\v} =
B_{\g_1}$ and $\nu$ is the representation of $\g_0$ on $\v$ defined as in (1.11).\vskip 1pc
By
Proposition 2.2 the Weyl algebra $W(\v) = S(\v)$ (as linear spaces) then inherits the
structure of a $\g_0$-module with respect to the action of $\tau(\nu(x))$ for $x\in \g_0$. But,
by Proposition 2.4, $\nu$ uniquely ``lifts" to a Lie algebra homomorphism $$\nu_*:\g_0\to
S^2(\v)\eqno (2.22)$$ such that, by (2.14), $$ad\,\nu_*(x) = \tau(\nu(x))\eqno (2.23)$$ for
any $x\in \g_0$. With respect to the bilinear forms  $B_{S(\v)}\vert S^2(\v)$ and $B_{\g_0}$ the
map $\nu_*$ admits a transpose $$\nu_*^t:S^2(\v)\to \g_0\eqno (2.24)$$ so that for any $x\in
\g_0$ and $w\in S^2(\v)$ one has $$(\nu_*(x),w) = (x,\nu_*^t(w))\eqno (2.25)$$ The spaces
$\g_0$ and $S^2(\v)$ are $\g_0$-modules respectively with respect to $ad$ and $ad\,\nu_*(x)$.
\vskip 1pc {\bf Proposition 2.7.} {\it Both $\nu_*$ and its transpose $\nu_*^t$ are
$\g_0$-maps.}\vskip 1pc {\bf Proof.} Clearly $\nu_*$ is a $\g_0$-map since $\nu_*$ is a Lie
algebra homomorphism. But then $\nu_*^t$ is a $\g_0$-map since $B_{\g_0}$ is, of course, 
$\g_0$-invariant and
$B_{S(\v)}$ is $\g_0$-invariant by Proposition 2.5. QED\vskip 1pc Let $U(\g_0)$ be the
universal enveloping algebra of $\g_0$. Since $\nu_*$ is a Lie algebra homomorphism, it extends
to an algebra homomorphism $$\nu_*:U(\g_0)\to W^{even}(\v)\eqno (2.26)$$ Let
$k=dim\,\g_0$ and let $x_i,\,i=1,\ldots,k$, be an orthonormal basis of $\g_0$ with respect to
$B_{\g_0}$. Consider the Casimir element $Cas_{\g_0} = \sum_{i=1}^k x_i^2$. Of course
$Cas_{\g_0}$ is in the center
$Z(\g_0)$ of $U(\g_0)$. We will be particularly interested in $\nu_*(Cas_{\g_0})$. By (2.20)
one has $$\nu_*(Cas_{\g_0}) \in S^4(\v) + \Bbb C \eqno (2.27)$$ 

2.4. We retain the assumptions and notation of \S 2.3 so that  $\g_0$ is an arbitrary (ordinary)
Lie algebra which admits a non-singular $ad$-invariant symmetric bilinear form
$B_{\g_0}$. Also $\v$ is an even dimensional
vector space with a non-singular alternating bilinear form $B_{\v}$ and $$\nu:\g_0\to
Lie\,Sp(\v)$$ is a symplectic representation (with respect to $B_{\v}$) of $\g_0$ on $\v$.
Recalling Remark 2.7 our main theorem will be a condition on $\nu$ which characterizes the case
when the assumptions arise from a Lie superalgebra of Riemannian type. 

Let
$\g = \g_0 \oplus \v$ and regard $\g$ as a $\Bbb Z_2$-graded vector space where $\g_0$ and
$\v$ are respectively the even and odd components of $\g$. Let $B_{\g}$ be the
non-singular supersymmetric bilinear form on $\g$ defined so that $\g_0$ and $\v$ are orthogonal
and
$B_{\g}\vert
\g_0 = B_{\g_0},\,\,B_{\g}\vert
\v = B_{\v}$. We will say that the pair $(\nu,B_{\g})$ is of super Lie type in case there
exists a Lie superalgebra structure on $\g$ such that (1) $\g_0$ is the even component, $\v=
\g_1$ is the odd component, $[x,y] = \nu(x)y$ for any $x\in \g_0,\,y\in \v$ and (2) $B_{\g_0}$
is $ad$-invariant. \vskip 1pc {\bf Theorem 2.8.} {\it Let $(\nu,B_{\g})$ be
arbitrary. Then $(\nu,B_{\g})$ is of super Lie type if and only if $\nu_*(Cas_{\g_0})$ is a 
constant in the Weyl algebra $W(\v)$. That is, $(\nu,B_{\g})$ is of super Lie type if and
only if, in the notation of (2.27), the component of
$\nu_*(Cas_{\g_0})$ in $S^4(\v)$ vanishes. Furthermore the  Lie superalgebra structure on
$\g$ is unique. In fact for any $y,y'\in \v$ one has $$[y,y'] = 2\nu_*^t(y\cdot y')\eqno
(2.28)$$ where we recall
$y\cdot y'\in S^2(\v)$ is the (commutative) product of $y$ and $y'$ in $S(\v)$ and 
$[y,y']\in \g_0$ is the supercommutator in $\g$ (and is not commutation in $W(\v)$). In addition
$\g$ is a Lie superalgebra of Riemannian type and $B_{\g}$ is a super-Riemannian bilinear form on
$\g$.}
\vskip 1pc {\bf Proof.} Assume first that $(\nu,B_{\g})$ is of super Lie type and let $[y,z]$
be the bracket in the Lie superalgebra $\g$. We will show (2.28) is necessarily satisfied. Let
$y,y'\in \v$ and let $x\in \g_0$. Then by (2.3), (2.11) and (2.23) one has
$$\eqalign{([y,y'],x)&= -(y',[x,y])\cr &= 2(y', \iota(y)\nu_*(x))\cr &= 2(y\cdot y',\nu_*(x))\cr
&= 2(\nu_*^t(y\cdot y'),x)\cr}\eqno (2.29)$$ But this establishes (2.28) since $x\in \g_0$
is arbitrary and $B_{\g_0}$ is non-singular. 

Now assume $(\nu,B_{\g})$ is arbitrary. We will define a bracket $[y,z]$ in $\g$ which will
clearly satisfy
$$[y,z] = - (-1)^{|y| |z|}[z,y]\eqno (2.30)$$  Let $[y,z]$ be the bracket in $\g_0$ if
$y,z\in\g_0$. Let 
$[y,z]= -[z,y] = \nu(y)z$ if $y\in \g_0$ and $z\in \v$. Let $[y,z] = 2\nu_*^t(y\cdot z)$ if
$y,z\in \v$. Now for any $x\in \g$ let $a(x)\in End\,\g$ be defined so that $a(x)y =
[x,y]$ for any $y\in\g$. We next show that $B_{\g}$ is invariant under $a(x)$ any $x\in \g$. That
is, for homogeneous $x,y,z\in \g$ one has $$([x,y],z) = - (-1)^{|x| |y|}(y,[x,z])\eqno (2.31)$$
If
$x\in \g_0$ then both sides vanish if $y$ and $z$ have unequal degrees. But one has (2.31) for
$y,z\in\g_0$ since $B_{\g_0}$ is $ad$-invariant and one has (2.30) for $y,z\in \v$ since
$\nu(x)\in Lie\,Sp(\v)$. Next assume $x\in \v$. Then both sides of (2.31) vanish if $y$ and
$z$ have the same degrees. By reversing sides it suffices only to consider only one of the
remaining 2 cases. Namely it suffices to consider the case where $y\in
\v$ and $z\in \v$. But then (2.31) follows from a reverse of the argument in (2.29). That is
going from bottom to the top. This establishes (2.31).

We now determine to what extent the super Jacobi identity $$[x,[y,z]] = [[x,y],z] + (-1)^{|x|
|y|} [y,[x,z]]\eqno (2.32)$$ is already satisfied. Of course (2.32) is satisfied in case
$x,y,z\in \g_0$. In case 2 of the 3 elements $x,y,z$ are in $\g_0$, then (2.32) is also satisfied
since $\nu$ is a representation. But (2.32) is also satisfied if only one of the elements
$x,y,z$ lies in $\g_0$. To prove this it suffices only to consider the case where $x\in \g_0$
and $y,z\in \v$. The fact that this case is sufficient uses (2.30) which holds for all
homogeneous pairs $x,y$ in $\g$.  But by Proposition 2.2, (2.23) and Proposition 2.7
$$\eqalign{[x,[y,z]] &= 2[x,\nu_*^t(y\cdot z)]\cr
&= 2\nu_*^t (ad\,\nu_*(x)(y\cdot z))\cr
&= 2\nu_*^t([x,y]\cdot z + y\cdot [x,z])\cr
&= [[x,y],z] + [y,[x,z]]\cr}$$ This establishes (2.32). The final case to be considered is
when the three elements are in $\v$. Recalling (B) in \S 1.1 we will have established that
$(\nu,B_{\g})$ is of super Lie type if and only if for all $y,y',y''\in \v$ one has $$
[y,[y',y'']] + [y',[y'',y]] + [y,[y',y'']] = 0\eqno (2.33)$$ Let $y,y',y''\in \v$. Let $w\in
S^2(\v)$ and $y,y',y''\in \v$. Then by Proposition 2.2 $$\eqalign{\iota(y')\iota(y)w^2&=
\iota(y')((\iota(y)w)\,w + w\,
\iota(y)w)\cr &= 2( \iota(y')\,\iota(y)w)\,w +\iota(y)w\, \iota(y')w + \iota(y')w\,
\iota(y)w\cr}\eqno (2.34)$$ But $\iota(u)\iota(v)w = (u\cdot v,w)$ by (2.3). Hence applying
$\iota(y'')$ to (2.34) yields
$$\iota(y'')\iota(y')\iota(y)w^2= 2((y\cdot y',w)\iota(y'')w + (y'\cdot y'',w)\iota(y)w +
(y''\cdot y, w)\iota(y')w)\eqno (2.35)$$ Now let $\{x_i\}, i=1,\ldots,k$, be a
$B_{\g_0}$-orthonormal basis of $\g_0$ so that $Cas_{\g_0} = \sum_{i=1}^k x_i^2$. But now if we
replace $w$ in (2.35) by $\nu_*(x_i)$ one has by (2.11) and (2.23)
$$\eqalign{\iota(y'')\iota(y')\iota(y)\nu_*(x_i)^2&= 2((y\cdot
y',\nu_*(x_i))\iota(y'')\nu_*(x_i)  + (y'\cdot
y'',\nu_*(x_i))\iota(y)\nu_*(x_i)\cr&\qquad\qquad + (y''\cdot y,
\nu_*(x_i))\iota(y')\nu_*(x_i))\cr &=  (\nu_*^t(y\cdot
y'),x_i)[y'',x_i]  + (\nu_*^t(y'\cdot
y''),x_i)[y,x_i]\cr&\qquad\qquad + (\nu_*^t(y''\cdot
y),x_i)[y',x_i]\cr &= {1\over 2}(([y,y'],x_i)[y'',x_i]  + ([y',y''],x_i)[y,x_i] +
([y'',y],x_i)[y',x_i])\cr }$$ But then summing over $i=1,\ldots,k$, yields
$$\iota(y'')\iota(y')\iota(y)\nu_*(Cas_{\g_0}) = {1\over 2}([y'',[y,y']] +[y,[y',y'']] +
[y',[y'',y]])\eqno (2.36)$$ But if $z\in S^4(\v)$ then for any
$y'''\in \v$, $(y''',\iota(y'')\iota(y')\iota(y)z) = (y\cdot y'\cdot y''\cdot y''',z)$ by (2.3).
Hence
$z=0$ if and only if $\iota(y'')\iota(y')\iota(y)z = 0$ for all $y,y',y''\in \v$ since
$B_{S(\v)}$ is non-singular on $S^4(\v)$. Thus one has the Jacobi identity (2.33) if and only
if the component of $\nu_*(Cas_{\g_0})$ in $S^4(\v)$ vanishes. The result then follows from
(2.27). QED \vskip 1pc Assuming $(\nu,B_{\g})$ is of super Lie type the scalar
$\nu_*(Cas_{\g_0})$ describes the action of $\nu_*(Cas_{\g_0})$ in the
representation (0.9) of $W(\v)$. We can compute this scalar purely in terms of the finite
dimensional representation $\nu$ of $\g_0$. Since
$Lie\,Sp(\v)$ is a simple Lie algebra, it has a unique
$ad$-invariant symmetric bilinear form up to scalar multiplication. However two such forms
present themselves to us. One,
$tr\,\alpha\beta$ for
$\alpha,\beta\in Lie\,Sp(\v)$, and the second is given by $B_{S_{\v}}\vert S^2(\v)$. The
following lemma computes the ratio. \vskip 1pc {\bf Lemma 2.9.} {\it Let $w,z\in S^2(\v)$. Then
$$(w,z) = {1\over 8}\,\,tr\,(ad\,w\,\,ad\,z)\vert \v\eqno (2.37)$$} \vskip 1pc {\bf Proof.} By
the simplicity of $Lie\,Sp(\v)$ it suffices to verify (2.37) for a single choice of $w$ and $z$
where $(w,z)\neq 0$. In fact we choose $w=z$ and $w = u\cdot v$ where $u,v\in\v$ are chosen so
that
$(u,v)= 1$. Since of course $(u,u) = (v,v) = 0$ it follows immediately from (2.2)
that $$(u\cdot v,u\cdot v) = 1\eqno (2.38)$$ On the other hand by (2.11) one readily has
$$\eqalign{ad\,u\,\, ad\,v(v) &= 2v\cr ad\,u\,\,ad\,v(u)&=-2u\cr}\eqno (2.39)$$ But
$ad\,u\,\, ad\,v$ clearly vanishes on the $B_{\v}$-orthocomplement of the plane spanned by $u$
and $v$. Hence $tr\,(ad\,u\,\, ad\,v)^2 = 8$. QED \vskip 1pc The representation $\nu$ of
$\g_0$ of course extends to a representation of $\nu$ of $U(\g)$ on $\v$. We can now prove
\vskip 1pc {\bf Theorem 2.10}. {\it Assume that $(\nu,B_{\g})$ is of super Lie type. Then 
$$\nu_*(Cas_{\g_0}) = {1\over 8}\,\,tr\,\nu(Cas_{\g_0})\eqno (2.40)$$}\vskip 1pc {\bf Proof.} By
Theorem 2.8, (2.15) and Lemma 2.9 one has, using notation in the proof of Theorem 2.8,
$$\eqalign{\nu_*(Cas_{\g_0}) &=
\chi(\nu_*(Cas_{\g_0}))\cr &= \sum_{i=1}^k (\nu_*(x_i),\nu_*(x_i))\cr
&= {1\over 8}\,\sum_{i=1}^k tr\, \nu(x_i)^2\cr &= {1\over 8}\,tr\,\nu(Cas_{\g_0})\cr}$$
QED\vskip 1.5pc


                          \centerline{\bf References}\vskip 1pc
\parindent=42pt

\item {[Con]} L. Conlon, A class of variationally complete
representations, {\it J. Diff. Geometry}, {\bf 7}(1972), 149-160
\item {[Kac]} V. Kac, Lie superalgebras, {\it Adv. in Math.}, {\bf
26}(1977), 8-96
\item {[Kos]} B. Kostant, A cubic Dirac operator and the emergence
of Euler number multiplets of representations for equal rank
subgroups, {\it Duke Math. Journal}, {\bf 100}(1999), 447-501
\item {[Parth]} R. Parthasarathy, Dirac operator and the discrete
series, {\it Ann. of Math.}, {\bf 96}\break (1972), 1-30
\smallskip
\parindent=30pt
\baselineskip=14pt
\vskip 1.9pc
\vbox to 60pt{\hbox{Bertram Kostant}
      \hbox{Dept. of Math.}
      \hbox{MIT}
      \hbox{Cambridge, MA 02139}}\vskip 1pc

      \noindent E-mail kostant@math.mit.edu
\end